\documentclass{elsart}
\usepackage{amsmath,amsfonts}

\journal{J.\ Math.\ Anal.\ Appl.}

\date{30 July 2002}

\newcommand{\ep}{\epsilon}
\newcommand{\supp}{\operatorname{supp}} 

\begin{document}

\begin{frontmatter}

\title{Singular measures and convolution operators}
\thanks{Research supported by grant BFM2000-0206-C04-03}

%
\author{J. M. Aldaz},
\ead{aldaz@dmc.unirioja.es}
\author{Juan L. Varona}
\ead{jvarona@dmc.unirioja.es}

\address{Departamento de Matem\'aticas y Computaci\'on,
Universidad de La Rioja,\\
26004 Logro\~no, Spain}


\begin{abstract}
We show that in the study of certain convolution operators,
functions can be replaced by measures without changing the size of the
constants appearing in weak type $(1,1)$ inequalities. As an application,
we prove that the best constants for the centered
Hardy-Littlewood maximal operator associated to parallelotopes do not
decrease with the dimension.
\end{abstract}


\end{frontmatter}

\section{Introduction}

The method of discretization for convolution
operators, due to M.~de Guzm\'an
(cf.~\cite{Guz}, Theorem~4.1.1), and further developed by
M.~T.~Men\'ar\-guez and F.~Soria (cf.~Theorem~1 of~\cite{MeSo})
consists in replacing functions by finite sums of Dirac
deltas in the study of the operator.
So far, the main applications of these theorems have been
related to the Hardy-Littlewood maximal function, and
more precisely, to the determination of bounds for the best
constants $c_d$ appearing in the weak type $(1,1)$ inequalities
(cf.~\cite{MeSo},~\cite{A1},~\cite{Mel1}, and \cite{Mel2} for the
one dimensional case, and for higher dimensions,
\cite{MeSo} and~\cite{A2}).
In this paper we complement de Guzm\'an's Theorem by proving
that one can
consider arbitrary measures instead of finite discrete measures, and the same
conclusions still hold (Theorem~\ref{th:equiv}). A special case of our
theorem (where the space is the real line and the convolution operator
is precisely the Hardy-Littlewood maximal function) appears
in~\cite{Mel2} (see Theorem~2).

Regarding upper bounds
for $c_d$, E.~M.~Stein and J.~Str\"omberg (see~\cite{StSt}) showed
that the constants grow at most like $O(d \log d)$ for arbitrary balls, and
like $O(d)$ in the case of euclidean balls. With respect to lower
bounds for the maximal function associated to cubes, it is shown
in~\cite{MeSo}, Theorem~6,
that $c_d \ge \big( \frac{1+2^{1/d}}{2} \big)^d$. These
bounds, which decrease with the dimension to $\sqrt{2}$,
where conjectured to be optimal in~\cite{Me}.
The ``optimality  part" of the conjecture was refuted
in~\cite{A2}, where
it was proved that $\liminf_d c_d \ge \frac{47 \sqrt{2}}{36}$.
It is an easy consequence of Theorem~\ref{th:equiv}
that the ``decreasing part'' of the conjecture
is also false:
For cubes the inequality $c_d\le c_{d+1}$ holds in every dimension $d$
(Theorem~\ref{th:ineq}). In dimensions $1$ and $2$ the stronger result
$c_1 < c_2$ is known, thanks to the recent determination by Antonios
D.~Melas of the exact value of $c_1$ as $\frac{11
+ \sqrt{61}}{12}$ (Corollary~1 of \cite{Mel2}). Since  $c_2 \ge
\sqrt\frac32 + \frac{3-\sqrt2}4$, by Proposition~1.4 of \cite{A2},
Melas's result entails that the first inequality is strict.

Finally, we note that the original question of Stein and Str\"omberg (see
also~\cite{BH}, Problem~7.74~c, proposed by A.~Carbery) as to whether
$\lim_d c_d < \infty$ or $\lim_d c_d = \infty$, remains open.

\section{Convolution operators and measures}

We shall state the main theorem of this note in terms of a
locally compact group~$X$.
Denote by $C(X)$ the family of all continuous
functions $g \colon X\to \mathbb{R}$, by $C_c(X)$ the continuous
functions with compact support, and by $\lambda$ the left Haar
measure on~$X$. If $X = \mathbb{R}^d$,
$\lambda^d$ will stand for the $d$-dimensional Lebesgue
measure. As usual, we shall write $dx$ instead of $d\lambda (x)$.
A finite real valued Borel measure $\mu$ on $X$ is \emph{Radon} if
$|\mu|$ is inner
regular with respect to the compact sets. It is well
known that if $X$ is a locally compact separable metric space, then
every finite Borel measure is automatically Radon.
Let $\mathcal{N}$ be a neighborhood base at $0$ such that each element
of $\mathcal{N}$ has compact closure, and let
$\{h_U : U \in \mathcal{N}\}$ be an approximate identity, i.e.,
a family of nonnegative Borel functions such that for every
$U \in \mathcal{N}$, $\supp h_U \subset U$ and $\|h_U\|_1 =1$.
Furthermore, since for every neighborhood $U$ of $0$ there
is a continuous function $g_U$ with values in $[0,1]$,
$g_U(0) = 1$, and $\supp g_U \subset U$, we may assume that each
function in the approximate identity is continuous (obtain $h_U$
by normalizing~$g_U$).
Let $\mu$ be a finite, nonnegative Radon measure on~$X$.
Recall that
\[
    h*f(x) = \int f(y^{-1}x) h(y) \, dy
    \text{\qquad and \qquad}
    \mu*f(x) = \int f(y^{-1}x) \, d\mu(y).
\]
Let $g\in C_c(X)$;
we shall utilize the following well known results:
$\mu*(h_U*g) = (\mu*h_U)*g$, and
$h_U*g \to g$ uniformly as $U \downarrow 0$.
The idea of the proof below consists simply in replacing the measure
$\mu$ with the continuous function $\mu*h_U$, using the fact that
$\| \mu*h_U \|_1 = \mu(X)$.

The $L_1$ norm refers always in this paper to Haar measure.

\begin{lem}
\label{lem:family}
Let $\{k_{\beta}\}$ be a family of
nonnegative lower semicontinuous real valued functions, defined
on~$X$. Set $k^*v := \sup_{\beta}|v*k_{\beta}|$, where $v$ is either a
function or a measure. Then, for every
finite real valued Radon measure $\mu$ on $X$,
and every $\alpha > 0$,
\[
    \lambda^d \{k^* \mu > \alpha\} \le
    \sup \big\{ \lambda^d \{k^* f > \alpha\} : \|f\|_1 = |\mu|(X) \big\}.
\]
The same result holds if $\{k_n\}$ is a sequence of
nonnegative real valued Borel functions.
\end{lem}

\begin{pf}
Consider first the case where $\{k_{\beta}\}$ is a family of
lower semicontinuous functions. We shall assume that
functions and measures are nonnegative.
There is no loss of generality in doing so since
$k^* f \le k^*|f|$ and
$k^*\mu \le k^*|\mu|$
always. Also,
by lower semicontinuity,
$\int k_{\beta} \, d\mu = \sup \{\int g_{\gamma,\beta} \, d\mu :
0\le g_{\gamma ,\beta} \le k_\beta$,
$g_{\gamma,\beta} \in C_c(X)\}$ (Corollary~7.13 of~\cite{Fo}).
It follows that for every $x$,
$\sup_\beta \mu *k_{\beta} (x) = \sup_{\gamma, \beta}
\{\mu *g_{\gamma,\beta} (x): 0\le g_{\gamma ,\beta} \le k_\beta$,
$g_{\gamma ,\beta} \in C_c(X)\}$. Therefore
we may assume that the family $\{k_{\beta}\}$
consists of nonnegative continuous functions with compact support.

Next, let $\{h_U: U \in \mathcal{N}\}$ be an approximate identity as above,
with each $h_U$ continuous, and
let $C \subset \{k^* \mu > \alpha\}$
be a compact set. It suffices to show that there exists
a function $f$ with $\| f \|_1 = \mu(X)$ and
$C \subset \{k^* f > \alpha\}$. We shall take $f$ to be $\mu*h_{U_0}$,
for a suitably chosen neighborhood $U_0$. Since
$\{k^* \mu > \alpha\} = \cup_\beta \{\mu *k_\beta > \alpha\}$ and each
$\mu *k_\beta $ is continuous, there exists a finite subcollection of
indices $\{\beta_1, \dots , \beta_\ell\}$ with
$C \subset \cup_1^\ell \{\mu *k_{\beta_i} > \alpha\}$, so the continuous
function $\max_{1\le i\le \ell} \mu *k_{\beta_i} $ attains a minimum
value $\alpha + a$ on $C$, with $a$ strictly positive.
Because $\mu$ is a finite measure and $h_U * k_{\beta_i}$
converges uniformly to $k_{\beta_i}$ as $U \to 0$,
$\mu * h_U * k_{\beta_i}$ also
converges uniformly to $\mu * k_{\beta_i}$. Hence, there exists an $U_0 \in
\mathcal{N}$ such that for every $V \subset U_0$, $V \in \mathcal{N}$, and
every $i\in \{1, \dots, \ell\}$,
$\| \mu * k_{\beta_i} - \mu * h_V * k_{\beta_i} \|_\infty < a /2$.
In particular, it follows that
\[
    C \subset \big\{ \max_{1 \le i \le \ell} \mu * h_{U_0} * k_{\beta_i}
      > \alpha \big\}
    \subset \big\{ k^*(\mu*h_{U_0})> \alpha \big\}.
\]

The case where $\{k_n\}$ is a sequence of
nonnegative bounded Borel functions, can be proven by reduction
to the previous one. Choose a finite Radon measure $\mu$
and fix $\alpha > 0$.
Given $\ep \in (0, 1)$, for every $n$ let
$g_n \ge k_n$ be a bounded, lower semicontinuous function with
\[
    \| g_n - k_n \|_1 < \frac{\ep^2}{2^{n+1}\mu(X)}
\]
(cf.~Proposition~7.14 of~\cite{Fo}).
Then, for any $ f \in L_1(\lambda)$, using the Fubini-Tonelli Theorem and
left invariance we have
\begin{gather*}
    \| g^*f - k^*f \|_1
\\
    = \Big\| \sup_n \int g_n (y^{-1}x) f(y) \,dy
       - \sup_n \int k_n (y^{-1}x) f(y) \,dy \Big\|_1
\\
    \le \sum_n \iint (g_n (y^{-1}x) - k_n (y^{-1}x)) |f(y)| \,dy\,dx
\\
    = \sum_n \int |f(y)| \int (g_n (y^{-1}x) - k_n (y^{-1}x)) \,dx\,dy
\\
    = \sum_n \|f\|_1 \|g_n - k_n\|_1
    < \|f\|_1 \ep^2 (\mu(X))^{-1}.
\end{gather*}
In particular, if $\| f \|_1 = \mu(X)$,
we have that
\[
    \| g^*f - k^*f \|_1 < \ep^2,
\]
from which
\[
    \lambda \{g^*f - k^*f \ge \ep\}
    \le \frac{\| g^* f - k^* f\|_1}{\ep} < \ep
\]
follows.
Now
$\{g^* f > \alpha + \ep \} \subset \{k^* f > \alpha\} \cup
\{g^* f - k^* f > \ep\}$, so
\begin{gather*}
    (\alpha + \ep) \lambda \{ k^*\mu > \alpha + \ep \}
    \le (\alpha + \ep ) \lambda \{ g^*\mu > \alpha + \ep \}
\\
    \le (\alpha + \ep) \sup \{ \lambda \{g^*f > \alpha + \ep\}
      : \|f\|_1 = \mu(X) \}
\\
    \le (\alpha + \ep) ( \sup \{\lambda \{k^*f > \alpha\}
      : \|f\|_1 = \mu(X)\} + \ep ),
\end{gather*}
and the result is obtained by letting $\ep \downarrow 0$.
\end{pf}

\begin{thm}
\label{th:equiv}
Let $\{k_{\beta}\}$ be a family of
nonnegative lower semicontinuous real valued functions, defined on
$X$, and let $c>0$ be a fixed constant. Then the
following are equivalent:

(i) For every function $f\in L_1 (\lambda )$,
and every $\alpha > 0$,
\[
    \alpha \lambda \{k^*f > \alpha\} \le c \|f\|_1.
\]

(ii) For every finite real valued Radon measure $\mu$ on $X$,
and every $\alpha > 0$,
\[
    \alpha \lambda \{k^*\mu > \alpha\} \le c |\mu|(X).
\]
The same result holds if $\{k_{n}\}$ is a sequence of
nonnegative real valued Borel functions.
\end{thm}

\begin{pf}
(i) is the special case of~(ii) where $d\mu(y) = f(y)\,dy$.
For the other direction, by Lemma~\ref{lem:family} and part~(i) we have
\[
    \alpha \lambda \{k^* \mu > \alpha\}
    \le \alpha \sup \{ \lambda \{k^*f > \alpha\} : \|f\|_1 = |\mu|(X) \}
    \le c |\mu|(X).
\]
\end{pf}

\begin{rem}
By the discretization theorem of M.~de Guzm\'an
(see~\cite{Guz}, Theorem~4.1.1), further refined by
M.~T.~Men\'arguez and F.~Soria (Theorem~1 of~\cite{MeSo}), in
$\mathbb{R}^d$ conditions~(i) and~(ii) of Theorem~\ref{th:equiv} are
both equivalent to

\textit{
(iii) For every finite collection $\{\delta_{x_1}, \dots , \delta_{x_N}\}$
of Dirac deltas on $X$, and every $\alpha > 0$,
\[
    \alpha \lambda \big\{ k^* \sum_1^N \delta_{x_i} > \alpha \big\} \le c N.
\]
}

From the viewpoint of obtaining lower
bounds, the usefulness of~(ii) is due to the fact that it allows to choose
among a wider class of potential examples than just finite sums of Dirac deltas.
Both~(ii) and~(iii) will be utilized in the sext section.
\end{rem}

\section{Behavior of constants for the Hardy-Littlewood maximal operator}
\label{sec:constants}

Let $B \subset \mathbb{R}^d$ be an open, bounded, convex set,
symmetric about zero.
We shall call $B$ a ball, since each norm on $\mathbb{R}^d$ yields sets of this
type, and each bounded $B$, convex and symmetric about zero, defines a norm.
The (centered) Hardy-Littlewood maximal operator associated to $B$
is defined for locally integrable
functions $f \colon \mathbb{R}^d \to \mathbb{R}$ as
\[
    M_{d,B} f(x) := \sup_{r > 0}
      \frac{\chi_{rB}}{r^d \lambda^d (B)}*|f| (x).
\]
We denote by $c_{d,B}$ the best constant in the weak type $(1,1)$
inequality
$\alpha \lambda^d \{ M_{d,B} f > \alpha \} \le c \|f\|_1$,
where $c$ is independent of
$f \in L^1(\mathbb{R}^n)$ and $\alpha > 0$. Let $s :=\{r_n\}_{-\infty}^\infty$
be a lacunary (bi)sequence (i.e., a sequence
that satisfies $r_{n+1}/r_n \ge c$ for some fixed constant
$c > 1$ and every $n \in \mathbb{Z}$). Then the associated maximal operator
is defined via
\[
    M_{s,d,B} f(x) := \sup_{n \in \mathbb{Z}}
    \frac{\chi_{r_n B}}{r_n^d \lambda^d(B)}*|f| (x).
\]
The arguments given below are applicable to both the maximal function
and to lacunary versions of it, so we shall not introduce
a different notation for the best constants in the lacunary case.
In particular, Lemma~\ref{lem:linear} and Theorem~\ref{th:ineq}
refer to all of these maximal operators, but only the usual maximal operator
shall be mentioned in the proofs.

Given a finite sum
$\mu = \sum_1^k \delta_{x_i}$ of Dirac deltas, where
the $x_i$'s need not be all different, let $\sharp (x + B)$ be the number
of point masses from $\mu$ contained in $x + B$.

\begin{lem}
\label{lem:linear}
Let $B$ be a ball in $\mathbb{R}^d$. Then for
every linear transformation $T \colon \mathbb{R}^d \to \mathbb{R}^d$ with
$\det T \ne 0$, $c_{d,B} = c_{d,T(B)}$.
\end{lem}

\begin{pf}
Given $\mu := \sum_1^k \delta_{x_i}$ and
$T \mu := \sum_1^k \delta_{T(x_i)}$, we have that
\[
    M_{d,B} \mu (x) := \sup_{r > 0} \frac{\sharp (x + rB)}{r^d\lambda^d (B)}
\]
and
\[
    M_{d,T(B)} T\mu (x) := \sup_{r > 0}
      \frac{\sharp (x + rT(B))}{r^d\lambda^d (T(B))}.
\]
Then
$x \in \{ M_{d,B} \mu > \alpha \}$
iff
$T(x)\in \{ M_{d,T(B)} T\mu > (\alpha / |\det T|)\}$.
Since
\[
    |\det T| \lambda^d \{ M_{d,B} \mu > \alpha \}
    = \lambda^d \{ M_{d,T(B)} T\mu > (\alpha / |\det T|) \},
\]
we have
\[
    \alpha \lambda^d \{ M_{d,B} \mu > \alpha \}
    = (\alpha / |\det T|) \lambda^d \{ M_{d,T(B)} T\mu
      > (\alpha / |\det T|) \},
\]
and the result follows.
\end{pf}

\begin{thm}
\label{th:ineq}
For each $d \in \mathbb{N} \setminus \{0\}$ let $B_d$ be
a $d$-dimensional parallelotope centered at zero. Then
$c_{d,B_d} \le c_{d+1,B_{d+1}}$ for both the maximal operator
and for lacunary operators.
\end{thm}

\begin{pf}
Since every such $B_d$ is the image under a nonsingular
linear transformation of the $d$-dimensional cube $Q_d$ centered at zero
with sides parallel to the axes and
volume $1$, we may assume that in fact $B_d = Q_d$.
With the convex bodies fixed, we will write
$c_d$ and $M_d$ rather than $c_{d,B_d}$ and $M_{d,B_d}$. Given $\alpha > 0$,
$\mu_d = \sum_1^k \delta_{x_i}$ on $\mathbb{R}^d$
and a constant $c > 0$ such that
$\alpha \lambda^d\{ M_d \mu_d > \alpha \} > c \mu_d (\mathbb{R}^d)$, we want to
find a measure $\mu_{d+1}$ on $\mathbb{R}^{d+1}$ such that
$\alpha \lambda ^{d+1}\{ M_{d+1} \mu_{d+1} > \alpha \} >
c \mu_{d+1} (\mathbb{R}^{d+1})$. This will imply that $c_d\le c_{d+1}$.
Let $L := (k/\alpha)^{1/d}$. Note that if $r \ge L$, then for every
$x\in \mathbb{R}^d$,
$\frac{\sharp (x + rQ_d) }{r^d} \le \alpha$. Choose $N \gg L$
such that $\alpha \frac{N-L}{N} \lambda^d \{ M_d \mu_d > \alpha \} > c k$,
and let $\mu_{d+1} := \mu_d \times \lambda_{[-N,N]}$, where
$\lambda_{[-N,N]}$ stands for the restriction of linear Lebesgue measure
to the interval $[-N, N]$. We claim that
$\{ M_d \mu_d > \alpha \} \times [-N+L, N-L] \subset
\{ M_{d+1} \mu_{d+1} > \alpha \}$. In order to establish the
claim, the following notation shall be used:
If $x = (x_1, \dots , x_d) \in \mathbb{R}^d$,
by $(x, x_{d+1})$ we denote the point
$(x_1, \dots , x_d, x_{d+1}) \in \mathbb{R}^{d+1}$. Now if
$x\in \{ M_d \mu_d > \alpha \}$, then there exists an $r(x)\in (0, L)$
such that $r(x)^{-d} \mu_d (x + r(x) Q_d) > \alpha$, so for every
$y\in [-N + L, N -L]$,
\begin{gather*}
    r(x)^{-d-1} \mu_{d+1} ((x, y) + r(x) Q_{d+1})
\\
    = r(x)^{-d-1} (\mu_d (x + r(x) Q_d)
      \times \lambda_{[-N,N]}( [y - \tfrac{r(x)}2, y + \tfrac{r(x)}2] )
\\
    = r(x)^{-d} \mu_d (x + r(x) Q_d)> \alpha,
\end{gather*}
as desired. But now
\begin{gather*}
    \alpha \lambda^{d+1} \{ M_{d+1} \mu_{d+1} > \alpha \}
    \ge 2 \alpha (N-L) \lambda^d \{ M_d \mu_d > \alpha \}
\\
    = 2 \alpha N \frac{N-L}{N} \lambda^d \{ M_d \mu_d > \alpha \}
       > 2 N c k = c \mu_{d+1} (\mathbb{R}^{d+1}).
\end{gather*}
\end{pf}

\begin{rem}
Recall from the Introduction that for the
$\ell_\infty$ balls (i.e., cubes with sides parallel to the axes)
$
c_1 < c_2.
$
Since the $\ell_1$ unit ball in dimension 2 is a square, it
follows from Lemma~\ref{lem:linear} that the best constant in dimension 2
is equal for the $\ell_1$ and the $\ell_\infty$ norms. It follows that
$c_1 < c_2$ in the $\ell_1$ case also. It would
be interesting to know whether or not the best constants associated
to the $\ell_p$ balls are all the same.
Note that  establishing bounds of the type
$a^{-1} c_{d,2} \le c_{d,p} \le a c_{d,2}$ (where the constant
$a \ge 1$ is independent of the dimension $d$ and $c_{d,p}$
denotes the best constant associated to the $\ell_p$ ball),
would show that the bounds
$O(d)$ (which hold for euclidean balls by
\cite{StSt}) extend to $\ell_p$ balls.
\end{rem}


\end{document}